# A Small Polyhedral Z-Acyclic 2-Complex in R4

## Authors

Frank H. Lutz and Günter M. Ziegler

## Description

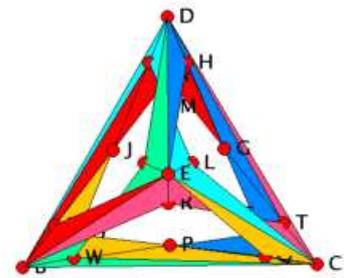

We present a small 4-dimensional polyhedral realization of a 2-dimensional Z-acyclic but non-contractible simplicial complex.

In November 2006, Lutz Hille (Hamburg) asked us for examples of Z-acyclic complexes that are not contractible, realized in low dimensions. This was motivated by cohomology computations for line bundles of toric varieties (cf. [1], [7], [8]), which can be reduced to cohomology computations for closed subcomplexes, as considered here. In this context, the question for the difference between homological triviality (acyclicity) and homotopical triviality (contractibility) arises naturally. Our model demonstrates explicitly that for 2-dimensional complexes piecewise linearly (PL) embedded in $\mathbf{R}^4$ (that is, for fans in $\mathbf{R}^5$, and thus for 5-dimensional toric varieties) homotopy conditions are stronger than homological ones. The respective dimensions hereby are smallest possible: Z-acyclic but non-contractible 2-dimensional complexes do not have PL embeddings in $\mathbf{R}^3$, whereas all 1-dimensional (finite) Z-acyclic complexes are contractible.

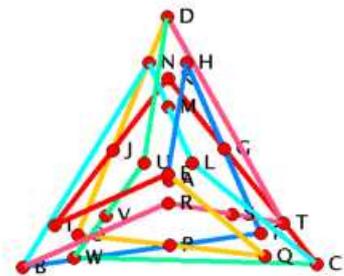

The most well known example of a Z-acyclic but non-contractible complex is due to Floyd and Richardson [4]. This 2-dimensional (cellular) complex K can be obtained by removing the open star of the point I from the quotient SO(3)/I of (an appropriate cell decomposition of) the group of 3-dimensional rotations SO(3) by the icosahedral group I. (The group I is isomorphic to the alternating group $A_5$ and acts freely on the complex K.) Alternatively, the complex K arises as the 2-skeleton of the Poincaré homology 3-sphere in its description by Weber and Seifert [16] as the spherical dodecahedral space. The complex K has the binary icosahedral group $I^*$ as its fundamental group.

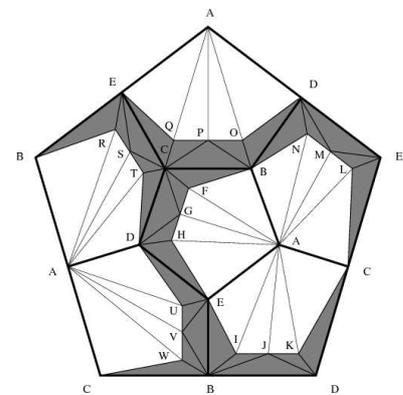

In fact, every homology 3-sphere $\Sigma^3$ different from $S^3$ can be used to yield a Z-acyclic but non-contractible 2-dimensional complex: One way to describe any such $\Sigma^3$ is as a 3-dimensional (polyhedral) ball with identifications on the boundary (cf. Seifert [15]). As in the Floyd-Richardson example, the identified boundary is a Z-acyclic but non-contractible 2-dimensional complex.

Newman [14] gave yet another construction of a Z-acyclic complex (non-homeomorphic to the complex K above, but also with fundamental group $I^*$) by identifying the boundaries of two circular discs. Furthermore, wedges as well as connected sums of Z-acyclic complexes are Z-acyclic again.

All (finite) 1-dimensional Z-acyclic complexes are trees and therefore are contractible. As planar graphs, they are geometrically realizable in $\mathbf{R}^2$.

The embedding respectively the PL embedding of 2-dimensional Z-acylic complexes is a more delicate matter. As pointed out by Zeeman [18], there are contractible 2-dimensional simplicial complexes that cannot be embedded (topologically) in any 3-manifold, for example, cones over non-planar graphs. Any (finite) 2-complex, however, is embeddable in some 4-manifold, and every contractible 2-complex can be embedded in $\mathbf{R}^4$; see Curtis [2]. (As remarked in [11], it is not known whether all contractible 2-dimensional complexes have a *PL* embedding in $\mathbf{R}^4$.) Such embeddings can be complicated (Curtis, cf. [5]): there is an embedding of a 2-dimensional contractible complex L in $S^4$ such that the complement $S^4 \setminus L$ is Z-acyclic (by Alexander duality) but not simply connected. For further embeddings with non-simply connected complement see Glaser [5], [6] and Neuzil [13].

It was shown by Kranjc [10; Cor. 2 and Rem. 2] that, in fact, arbitrary (finite) Z-acylic 2-dimensional complexes always are topologically embeddable in $\mathbf{R}^4$. A class of polyhedral homology 3-spheres in $\mathbf{R}^4$ with non-trivial fundamental group was given by Curtis and Wilder [3]: These give examples of PL embeddings of non-contractible Z-acylic 2-dimensional complexes in $\mathbf{R}^4$. In the following, we present a rather small model.

**Theorem:** There is a polyhedral realization (with small integer coordinates) in $\mathbf{R}^4$ of a 2-dimensional Z-acyclic but non-contractible simplicial complex with 23 vertices.

*Proof:* We appropriately subdivide the 2-skeleton of the spherical dodecahedral space (see above and cf. [12]) to obtain a Z-acyclic but non-contractible 2-dimensional simplicial complex with 23 vertices, A, B,...,W; see the bottom image for the respective triangulation. We then remove the star of the vertex A from the triangulation to obtain the displayed shaded complex. A 3-dimensional geometric realization of this shaded complex is given in the top applet. Finally, the star of vertex A is added by placing this vertex in $\mathbf{R}^4$ at position (0,0,0,1). To be more precise: The vertices B, C, D, and E of the shaded complex are connected to each other and therefore form the complete graph $K_4$ which is the 1-skeleton of the tetrahedron. Each of the six original pentagons of the 2-skeleton of the spherical dodecahedral space is subdivided (in a symmetric way) by inserting three additional vertices each. For example, the pentagon A-B-C-D-E is subdivided by adding the vertices F, G, and H. In each of the six pentagons there is a strip of five shaded triangles which is glued to the tetrahedral graph formed by the vertices B, C, D, and E. Each of the strips is colored in a different color, and it can be inspected in the applet how the strips wind arround the tetrahedral graph. The coordinates of the realization are chosen such that the rotation group of the tetrahedron acts on the model. All the coordinate entries that occur in the model are integers between -4 and +4.

The 3-dimensional realization of the shaded complex cannot be extended to a 3-dimensional realization of the full complex since the link of vertex A is linked in the realization (see the second applet).

In fact, 2-dimensional Z-acyclic but non-contractible simplicial complexes never have PL embeddings in $\mathbf{R}^3$: Suppose, there is such a PL embedding. Any regular neighborhood (cf. [17]) of the embedded complex yields a Z-acyclic 3-manifold with $\mathbf{S}^2$ as boundary (cf. [9, p. 114]). However, by the Schoenflies theorem, $\mathbf{S}^2$ bounds a 3-ball with trivial fundamental group, contradiction.

Acknowledgement: The authors were supported by the DFG Research Group "Polyhedral Surfaces", Berlin.

**Keywords**      Z-acyclic and contractible complexes; polyhedral realizations
**MSC-2000 Classification**    52B99 (57Q15, 51M20)

**Authors' Addresses**

Frank H. Lutz

    TU Berlin
    Institute of Mathematics, MA 3-2
    D-10623 Berlin
    lutz@math.tu-berlin.de
    http://www.math.tu-berlin.de/~lutz

Günter M. Ziegler

    TU Berlin
    Institute of Mathematics, MA 6-2
    D-10623 Berlin
    ziegler@math.tu-berlin.de
    http://www.math.tu-berlin.de/~ziegler